\documentstyle[bezier,12pt,leqno]{article}
\textheight9in \def\DATE{\today}
\def\cal{\mathcal}

\newtheorem{theorem}{Theorem}

\newtheorem{proposition}[theorem]{Proposition}
\newtheorem{remark}[theorem]{Remark}

\newtheorem{situation}[theorem]{Situation}
\newtheorem{problem}[theorem]{Problem}

\catcode`\@=11

\def\ps@myheadings{\let\@mkboth\@gobbletwo
\def\@oddhead{\ifnum\count0=1 \hfill\else
\rightmark \hfil \rm\thepage\fi}%
\def\@oddfoot{\ifnum\count0=1 \hfill \rm 1 \hfill \else
\hfill\fi}
\def\@evenhead%
{\rm\leftmark\hfil\rm\thepage}%
\def\@evenfoot{}\def\sectionmark##1{}
\def\subsectionmark##1{}}

\def\@begintheorem#1#2{\it \trivlist \item[\hskip
 \labelsep{\bf #1\ #2.}]}
\def\@opargbegintheorem#1#2#3{\it \trivlist\item[\hskip%
 \labelsep{\bf #1\ #2.\ (#3)}]}
\def\@endtheorem{\endtrivlist}

\def\@listI{\leftmargin\leftmargini \parsep 1pt plus 2.5pt
 minus 1pt\topsep 5pt plus 2pt minus 3pt%
 \itemsep 0pt plus 2.5pt minus 1pt}
\let\@listi\@listI
\@listi

\def\@sect#1#2#3#4#5#6[#7]#8{\ifnum #2>\c@secnumdepth%
 \def \@svsec {}\else \refstepcounter {#1}\edef \@svsec%
 {\csname the#1\endcsname. \hskip .1em }\fi \@tempskipa%
 #5\relax \ifdim \@tempskipa >\z@ \begingroup #6\relax%
 \@hangfrom {\hskip #3\relax \@svsec }{\interlinepenalty%
 \@M #8.\par }\endgroup \csname #1mark\endcsname {#7}%
 \addcontentsline {toc}{#1}{\ifnum #2>\c@secnumdepth%
 \else \protect \numberline {\csname the#1\endcsname. }%
 \fi #7}\else \def \@svsechd {#6\hskip #3\@svsec #8.%
 \csname #1mark\endcsname {#7}\addcontentsline {toc}{#1}%
 {\ifnum #2>\c@secnumdepth \else \protect \numberline%
 {\csname the#1\endcsname. }\fi #7}}\fi \@xsect {#5}}

\def\section{\@startsection {section}{1}{\z@ }%
 {-3.5ex plus -1ex minus -.2ex}{2.3ex plus .2ex}{\bf }}

\def\thebibliography#1{%
 \section *{References.\@mkboth {REFERENCES}{REFERENCES}}%
 \list {[\arabic {enumi}]}{\settowidth \labelwidth {[#1]}%
 \leftmargin \labelwidth \advance \leftmargin \labelsep %
 \usecounter {enumi}} \def \newblock %
 {\hskip .11em plus .33em minus -.07em} \sloppy \clubpenalty 4000%
 \widowpenalty 4000 \sfcode`\.=1000\relax}

\def\@maketitle{%
 \newpage \null \vskip 2em
 \begin{center}
{\Large\sf \@title \par }
 \vskip 1.5em
 {\large \lineskip .5em
 \begin {tabular}[t]{c}\@author
 \end{tabular}\par}
 \end{center}
  \vskip .8em}

\def\abstract{%
\if@twocolumn \section *{Abstract}
 \else \small\quotation\noindent{\bf Abstract.}\fi}

\catcode`\@=13

\def\qbezier{\bezier{100}}
\def\Tc{{T^c}} \def\Trp{{\tt P}} \def\Trq{{\tt Q}}
\def\p{\mbox{$p \hskip -2.15mm p \hskip - 2.15mm p$}}    \def\Trq{{\tt Q}}
\def\q{\mbox{$q \hskip -2.15mm q \hskip - 2.15mm q$}}    \def\Vert{{\it Vert}}
\def\ot{\otimes}  
\def\pa{\partial} \def\id{1\!\!1} \def\ainfty{{A_\infty}}
\def\phi{\varphi} \def\pkernels{p-kernels} \def\qkernels{q-kernels}
\def\otexp#1#2{#1^{\otimes #2}} 
\def\pkernel{p-kernel} \def\qkernel{q-kernel}
\def\orada#1#2{{#1\otimes\cdots\otimes#2}}
\def\bfmu{\mbox{$\mu \hskip -2.4mm \mu \hskip - 2.4mm \mu$}}
\def\bfxi{\mbox{$\xi \hskip -2.1mm \xi \hskip - 2.1mm \xi$}}
\def\bfeta{\mbox{$\eta \hskip -2.1mm \eta \hskip - 2.1mm \eta$}}
\def\bftheta{\mbox{$\theta \hskip -2.1mm \theta \hskip - 2.1mm \theta$}}
\def\bfnu{\mbox{$\nu \hskip -2.4mm \nu \hskip - 2.4mm \nu$}}
\def\bfphi{\mbox{$\phi \hskip -2.57mm \phi \hskip - 2.57mm \phi$}}
\def\bfpsi{\mbox{$\psi \hskip -2.73mm \psi \hskip - 2.73mm \psi$}}
\def\bfH{\mbox{$H \hskip -3.65mm H \hskip - 3.65mm H$}}
\def\bfC{\mbox{$C \hskip -3.4mm C \hskip - 3.4mm C$}}
\def\bfL{\mbox{$L \hskip -3mm L \hskip - 3mm L$}}
\def\susp{\uparrow\!}             \def\desusp{\downarrow\!}
       \def\Rada#1#2#3{#1_{#2},\dots,#1_{#3}}
\def\squeezedcdots{\hskip -.5mm \cdot \hskip -.5mm  \cdot \hskip -.5mm \cdot}
\def\Pin{{\cal P}_{\it in}} \def\Pout{{\cal P}_{\hskip -1mm \it out}}
\def\Rin{{\cal R}_{\it in}} \def\Rout{{\cal R}_{\hskip -1mm \it out}}
\def\rhoin{{\rho_{\it in}}} \def\rhoout{{\rho_{\it out}}}
\def\setinf{{\cal A}_\infty} \def\Trans#1#2#3{{\sf Tr}_{#1,#2,#3}}
\def\Hom{{\it Hom}}

\newcommand{\square}[8]{
\setlength{\unitlength}{.7cm}
\begin{picture}(5,3.6)
\thicklines

\put(0,3){\makebox(0,0){$#1$}}
\put(5,3){\makebox(0,0){$#2$}}
\put(0,0){\makebox(0,0){$#3$}}
\put(5,0){\makebox(0,0){$#4$}}

\put(-.5,1.5){\makebox(0,0)[r]{\scriptsize $#6$}}
\put(5.5,1.5){\makebox(0,0)[l]{\scriptsize $#7$}}
\put(2.5,0.5){\makebox(0,0)[b]{\scriptsize $#8$}}
\put(2.5,3.5){\makebox(0,0)[b]{\scriptsize $#5$}}

\put(1,0){\vector(1,0){3}}
\put(1,3){\vector(1,0){3}}
\put(0,2.5){\vector(0,-1){2}}
\put(5,2.5){\vector(0,-1){2}}
\end{picture}
}

\def\krouzekstopka{
\unitlength .9mm
\linethickness{0.4pt}
\ifx\plotpoint\undefined\newsavebox{\plotpoint}\fi 
\begin{picture}(0,3)
\put(0,4.23){\line(0,-1){2}}
\put(0,0){\line(0,-1){2}}
\put(0,1){\makebox(0,0)[cc]{\Large $\circ$}}
\end{picture}
}

\def\teckastopka{
\unitlength .9mm
\linethickness{0.4pt}
\ifx\plotpoint\undefined\newsavebox{\plotpoint}\fi 
\begin{picture}(0,3)
\put(0,4.23){\line(0,-1){2}}
\put(0,0){\line(0,-1){2}}
\put(0,1){\makebox(0,0)[cc]{\Large $\bullet$}}
\end{picture}
}
\def\qed{\hspace*{\fill}
\mbox{\hphantom{mm}\rule{0.25cm}{0.25cm}}\\}

\title{Transferring $A_\infty$ (strongly homotopy associative)
  structures}

\author{Martin Markl%
\thanks{The author was supported by Grant GA \v CR 201/02/1390.}}

\begin{document}

\maketitle

\begin{abstract}
  The aim of this simple-minded ``applied'' note is to give explicit
  formulas for transfers of $\ainfty$-structures and related maps and
  homotopies in the most easy situation in which these transfers
  exist.  The existence of these transfers follows, in characteristic
  zero, from a general theory developed by the author
  in~\cite{markl:ha}. The easier half of our formulas was already
  known to Kontsevich-Soibelman and
  Merkulov~\cite{kontsevich-soibelman:00,merkulov:98} who derived
  them, without explicit signs, under slightly stronger assumptions
  than those made in this note.
\end{abstract}

{\small
\vskip 3mm
\noindent 
{\bf \hskip 10mm  Table of content:} \ref{sec:1}.  
                     Introduction and results
                     -- page~\pageref{sec:1}
\hfill\break\noindent 
\hphantom{{\bf \hskip 10mm Table of content:\hskip .5mm}}  \ref{sec:2}.
                     Conventions
                     -- page~\pageref{sec:2}
\hfill\break\noindent 
\hphantom{{\bf \hskip 10mm Table of content:\hskip .5mm}} \ref{sec:3}.
                     Inductive formulas
                     -- page~\pageref{sec:3}
\hfill\break\noindent 
\hphantom{{\bf \hskip 10mm Table of content:\hskip .5mm}}  \ref{sec:4}.
                     Non-inductive formulas
                     -- page~\pageref{sec:4}
\hfill\break\noindent 
\hphantom{{\bf \hskip 10mm Table of content:\hskip .5mm}}  \ref{sec:5}. 
                     Why do the transfers exist
                     -- page~\pageref{sec:5}
\hfill\break\noindent 
\hphantom{{\bf \hskip 10mm Table of content:\hskip .5mm}}  \ref{sec:6}. 
                     Some other properties of the transfer
                     -- page~\pageref{sec:6}
\hfill\break\noindent 
\hphantom{{\bf \hskip 10mm Table of content:\hskip .5mm}}  \ref{sec:7}. 
                     Two observations
                     -- page~\pageref{sec:7}

\vskip 3mm

\noindent 
{\bf \hskip 10mm Classification:} 18D10; 55S99 
\hfill\break\vskip -5mm\noindent
{\bf \hskip 10mm Keywords:} $\ainfty$-algebra, transfer
}

\begin{center}
  -- -- -- -- -- 
\end{center}

\bibliographystyle{plain}

\section{Introduction and results}
\label{sec:1}

We will work in the category of (left) modules over an arbitrary
commutative unital ring~$R$. Therefore, by a chain complex we will
understand a chain complex of $R$-modules, by a linear map an
$R$-linear map, etc. In particular, results of this paper apply to the
category of abelian groups and to the category of vector spaces over a
field of arbitrary characteristic. Let us consider the following
situation.

\begin{situation}
  \label{sit}
  We are given chain complexes $(V,\pa_V)$, $(W,\pa_W)$ and chain maps
  $f : (V,\pa_V) \to (W,\pa_W)$, $g : (W,\pa_W) \to (V,\pa_V)$ such
  that the composition $gf$ is chain homotopic to the identity $\id_V:
  V\to V$, via a chain-homotopy $h$.
\end{situation}

A compact way to express Situation~\ref{sit} is to say that $g :
(W,\pa_W) \to (V,\pa_V)$ is a left chain-homotopy inverse of $f :
(V,\pa_V) \to (W,\pa_W)$.  Our assumptions are in particular satisfied
when the complexes $(V,\pa_V)$ and $(W,\pa_W)$ are chain homotopy
equivalent. In this note we address the following

\begin{problem}
\label{problem}
Suppose we are given an $\ainfty$-structure $\bfmu = (\mu_2,\mu_3, \ldots)$
on $(V,\pa_V)$. In Situation~\ref{sit}, give explicit formulas for
the following objects:
\begin{itemize}
\item [(i)]
an $\ainfty$-structure $\bfnu = (\nu_2,\nu_3, \ldots)$ on $(W,\pa_W)$,
\item [(ii)]
an $\ainfty$-map $\bfphi = (\phi_1,\phi_2,\ldots) : (V,\pa,\mu_2,\mu_3,
\ldots) \to (W,\pa,\nu_2,\nu_3, \ldots)$,
\item [(iii)]
an $\ainfty$-map $\bfpsi = (\psi_1,\psi_2,\ldots) : (W,\pa,\nu_2,\nu_3,
\ldots) \to (V,\pa,\mu_2,\mu_3, \ldots)$, and
\item [(iv)]
an $\ainfty$-homotopy $\bfH = (H_1,H_2,\ldots)$ between
$\bfpsi\bfphi$ and $\id_V$
\end{itemize}
such that $\bfphi$ extends $f$, $\bfpsi$ extends $g$ and $\bfH$
extends $h$ or, expressed more 
formally, $\phi_1 = f$, $\psi_1 = g$ and $H_1 = h$.
\end{problem}

Our strategy will be to construct suitable degree $n-2$ maps $\{\p_n :
\otexp Vn \to V\}_{n \geq 2}$ (the {\em \pkernels}) and suitable
degree $n-1$ maps $\{\q_n : \otexp Vn \to V\}_{n \geq 1}$ (the {\em
\qkernels}) such that $\nu_n$, $\phi_n$, $\psi_n$ and $H_n$ defined by
the following Anzatz:
\begin{equation}
\label{eq:1}
\nu_n : = f \circ \p_n \circ \otexp g n,\ \phi_n : = f \circ \q_n,\
\psi_n := h \circ \p_n \circ \otexp g n \mbox { and } H_n := h \circ
\q_n
\end{equation}
answer Problem~\ref{problem}.  We give both inductive
(formulas~(\ref{Dny_otevrenych_dveri})
and~(\ref{Eli_toto_utery_priletela}) in
Section~\ref{sec:inductive-formulas}) and non-inductive
(Propositions~\ref{pred_odletem_do_Izraele} 
and~\ref{po_Xmas_party_na_MFF_mne_boli_siska} of
Section~\ref{sec:non-induct-form}) formulas for the kernels.

\begin{remark}
  {\rm We already mentioned in the Abstract that the formulas for $\nu_n$
    and $\psi_n$ were given, without explicit signs,
    in~\cite{kontsevich-soibelman:00} (non-inductive formulas) and also
    in~\cite{merkulov:98} (inductive formulas).  Kontsevich and
    Soibelman~\cite{kontsevich-soibelman:00} assumed (in our notation)
    that $(W,\pa_W)$ was a subcomplex of $(V,\pa_V)$, $f : (V,\pa_V)
    \to (W,\pa_W)$ a projection, $g :(W,\pa_W) \hookrightarrow
    (V,\pa_V)$ the inclusion and, of course, that $gf$ was chain
    homotopic to the identity $\id_V$. Merkulov~\cite{merkulov:98}
    made similar assumptions and he moreover assumed that
    $(V,\pa_V,\bfmu)$ was an ordinary dg-associative algebra, that is,
    $\mu_n = 0$ for $n \geq 3$.  Our formulas for $\phi_n$ and $H_n$
    are, to our best knowledge, new ones.
    A~surprising interpretation of the \pkernel\ in terms of homotopy
    operads is suggested by~\cite{laan:03}.%
}
\end{remark}

\begin{remark}
\label{kourim_dymku}
{\rm
Given the data spelled out in Situation~\ref{sit}, 
suppose that these data satisfy,
furthermore, the {\em side conditions\/} (annihilation requirements):
 \[
f\circ h = 0,\ h\circ g = 0\ \mbox { and }\ h\circ h = 0.
\]
By means of a suitable modification of the data if need be,
these conditions can always be arranged for, see for instance the
discussion in~\cite{markl:ip}.
Then an application of the Coalgebra Perturbation
Lemma of Huebschmann and
Kadeishvili~\cite[$2.1_*$]{huebschmann-kadeishvili:MZ91}
yields a solution of the transfer problem spelled out as Problem~\ref{problem}.
A closer look reveals that, indeed,
kernels are lurking behind
the formulas of~\cite{huebschmann-kadeishvili:MZ91}
as well but in that paper there was no need to
introduce explicitly  objects equivalent to our kernels.
The r\^ole of the \pkernel\ is played by the
summation $\sum_{n \geq 0} (\tilde h \circ \delta_\mu)^n$ and the
\qkernel\ is represented by $\sum_{n \geq 0} (\delta_\mu \circ \tilde
h)^n$, where $\delta_\mu$ is the square-zero coderivation of
$\Tc(\desusp V)$ corresponding to the $A_\infty$-structure $\bfmu$ and
$\tilde h$ is the extension of $h$ as a coderivation homotopy,
see~\cite[Perturbation Lemma~1.1]{huebschmann-kadeishvili:MZ91}. It
can be shown that under the circumstances of
that paper, that is, when the side conditions are met,
these kernels coincide with the kernels
to be used in the present
paper. Without the side conditions, the formulas
of~\cite{huebschmann-kadeishvili:MZ91}
do not apply directly. For an historical overview of the problem of
transferring $A_\infty$ structures and other useful information
see~\cite{huebschmann:08}. 
}
\end{remark}
                                 
\begin{center}
  -- -- -- -- --
\end{center}

This work was stimulated by E.~Getzler who indicated that there might
be some need for explicit transfers.  The $\ainfty$-case discussed
here in fact turned out to be more elementary than we expected, which
we attribute to the existence of a canonical non-$\Sigma$
polarization~\cite[Remark~25]{markl:ho}. In a forthcoming paper we
intend to make similar constructions for strongly homotopy Lie
algebras.

\vskip 2mm

\noindent 
{\bf Acknowledgment:} We are indebted to our wife Kv\v etoslava for
sketching out the carp's head on page~\pageref{fish}, and to Johannes
Huebschmann for clarifying the relation between the present paper
and~\cite{huebschmann-kadeishvili:MZ91}.

\section{Conventions}
\label{sec:conventions}
\label{sec:2}

In this unbelievably boring section we set up sign conventions
used in this note. The signs in the axioms of $\ainfty$-algebras and
related objects are unique up to an action of the infinite product
$\mbox {\Large $\times$}_1^\infty C_2$ of the cyclic group $C_2 =
\{-1,1\}$. For example, $(\epsilon_2,\epsilon_3,\ldots) \in C_2 \times
C_2 \times \cdots$ acts on the signs in Axiom~(\ref{X}) below by
\[
(\mu_2,\mu_3,\ldots) \longmapsto (\epsilon_2 \mu_2, \epsilon_3
\mu_3,\ldots).
\]

The sign convention used here is compatible with the one
of~\cite{markl:JPAA92}. It differs from the original one of Jim
Stasheff~\cite{stasheff:TAMS63} by the action, in Axiom~(\ref{X}), of
$(\epsilon_2,\epsilon_3,\ldots)$ with $\epsilon_n = (-1)^{n(n-1)/2} =
\hbox{$\otexp \susp n \circ \otexp \desusp n$}$, where $\susp$
(resp.~$\desusp$\hskip .5mm) denotes the suspension
(resp.~desuspension) operator.

We are going to recall axioms for an $\ainfty$-structure $\bfmu =
(\mu_2,\mu_3,\cdots)$ on $(V,\pa_V)$ (Axiom~(\ref{X})), an
$\ainfty$-structure $(\nu_2,\nu_3,\ldots)$ on $(W,\pa_W)$
(Axiom~(\ref{Z})), for an $\ainfty$ map $\bfphi : (V,\pa_V,\bfmu) \to
(W,\pa_W,\bfnu)$ (Axiom~(\ref{Y})), for an $\ainfty$ map $\bfpsi :
(W,\pa_W,\bfnu) \to (V,\pa_V,\bfmu)$ (Axiom~(\ref{A})) and for an
$\ainfty$-homotopy $\bfH$ between the composition $\bfpsi\bfphi$ and
$\id_V$ (Axiom~(\ref{B})). In Axioms~(\ref{X}) and~(\ref{Z}), $\mu_n : \otexp
Vn \to V$ and $\nu_n : \otexp Wn \to W$ are $n$-multilinear degree
$n-2$ maps, in Axioms~(\ref{Y}) and~(\ref{A}), $\phi_n : \otexp Vn \to W$ and
$\psi_n : \otexp Wn \to V$ are $n$-multilinear maps of degree $n-1$,
and finally in Axiom~(\ref{B}), $H_n : \otexp Vn \to V$ is an $n$-multilinear
degree $n$ map. Here are the axioms in their full glory:

\begin{eqnarray}
\label{X}
\delta (\mu_n)
\hskip -2mm &:=& \hskip -2mm
\sum_A
(-1)^{i(l+1) +n} \mu_k(\id^{\ot i-1} \ot \mu_l \ot \id^{\ot k-i}),\ 
n \geq 2,
\\
\label{Z}
\delta (\nu_n)
\hskip -2mm &:=& \hskip -2mm
\sum_A
(-1)^{i(l+1) +n} \nu_k(\id^{\ot i-1} \ot \nu_l \ot \id^{\ot k-i}),\
n \geq 2,
\\
\label{Y}
\delta (\phi_n) \hskip -2mm &:=& \hskip -2mm
-\sum_B (-1)^{\vartheta(r_1,\ldots,r_k)}
\nu_k(\orada{\phi_{r_1}}{\phi_{r_k}})+
\\
\nonumber 
&&\hskip 1cm - \sum_A
(-1)^{i(l+1) +n}
\phi_k(\id^{\ot i-1} \ot \mu_l \ot \id^{\ot k-i}),\ n \geq 1,
\\
\label{A}
\delta (\psi_n) \hskip -2mm &:=& \hskip -2mm
-\sum_B (-1)^{\vartheta(r_1,\ldots,r_k)}
\mu_k(\orada{\psi_{r_1}}{\psi_{r_k}})+
\\
\nonumber 
&&\hskip 1cm - \sum_A
(-1)^{i(l+1) +n}
\psi_k(\id^{\ot i-1} \ot \nu_l \ot \id^{\ot k-i}),\ n \geq 1, \mbox { and}
\\
\label{B}
\delta (H_n) \hskip -2mm &:=& \hskip -2mm
- \sum_C
(-1)^{n + r_i + 
\vartheta(r_1,\ldots,r_i)
}
\mu_k((\bfpsi\bfphi)_{r_1}\! \ot
\squeezedcdots \! \ot \!
(\bfpsi\bfphi)_{r_{i-1}}  \!\ot \! H_{r_i} \! \ot \! \otexp {\id}{k-i})
\\
\nonumber 
&& \hskip .7cm
+ \sum_A (-1)^{n + i(l+1)} 
H_k (\otexp {\id}{i-1}\! \ot\! \mu_l\! \ot\! \otexp {\id}{k-i}) + (\bfpsi
\bfphi)_n\! -\! (\id_V)_n,\
n \geq 1.
\end{eqnarray}
In the above display, 
\begin{eqnarray*}
A \hskip -2mm &:=& \hskip -2mm \{k,l\ |\ k+l = n+1,\ k,l \geq 2,\ 
1 \leq i \leq k\},
\\
B \hskip -2mm &:=& \hskip -2mm \{k, \Rada r1k\ |\ 2 \leq k \leq n,\ 
r_1,\ldots, r_k \geq 1,\ r_1 + \cdots + r_k = n\},
\\
C \hskip -2mm &:=& \hskip -2mm \{k,i,\Rada r1i\ |\   
       2 \leq k \leq n,\ 1 \leq i \leq k,\ r_1,\ldots, r_i \geq 1,\
r_1+ \cdots + r_i + k - i = n\},
\end{eqnarray*}
and, for integers $u_1,\ldots,u_s$, we denoted
\[
\vartheta(u_1,\ldots,u_s) := \sum_{1 \leq \alpha < \beta \leq s}
u_\alpha(u_\beta + 1).
\]
The symbols $\delta$ in the left hand sides denote the induced
differentials in the corresponding complex of multilinear maps.  To
interpret the above axioms in terms of elements, one must of course use
the Koszul sign convention. For example, Axiom~(\ref{X}) evaluated at elements
$\Rada v1n \in V$, reads
\begin{eqnarray*}
\lefteqn{\hskip -1cm
\pa_V \mu_n(\Rada v1n) - 
\sum_{1 \leq i \leq n}
(-1)^{n+ |v_1| + \cdots + |v_{i-1}|}\mu_n(\Rada v1{i-1},\pa_V (v_i),\Rada
v{i+1}n) := 
}
\\
\hskip 1cm &&:=
\sum_A
(-1)^{i(l+1) +n + l(|v_1| + \cdots + |v_{i-1}|)} 
\mu_k(\Rada v1{i-1},\mu_l(\Rada vi{i+l-1}),\Rada v{i+l}n),
\end{eqnarray*}
which is~\cite[Equation~(1)]{markl:JPAA92}.
If $\Tc(-)$ denotes the tensor coalgebra functor, then
\begin{itemize}
\item
$\bfmu$ is the same as a degree~$-1$ square-zero coderivation
$\delta_\mu$ of $\Tc(\desusp V)$ whose linear part is $\pa_V$, 
\item
$\bfnu$
is the same as a degree~$-1$ square-zero coderivation $\delta_\nu$ of
$\Tc(\desusp W)$ with linear part~$\pa_W$, 
\item
$\bfphi$ is
the same as a dg-algebra homomorphism $F : (\Tc( \desusp V),\delta_\mu) \to
(\Tc(\desusp W),\delta_\nu)$, 
\item
$\bfpsi$ is the same as a dg-algebra
homomorphism $G : (\Tc(\desusp W),\delta_\nu) \to (\Tc(\desusp V),\delta_\mu)$,
and 
\item
$\bfH$ is the same as a coderivation homotopy between $GF$ and the
identity map of~$\Tc(\desusp V)$.
\end{itemize}

\section{Inductive formulas}
\label{sec:inductive-formulas}
\label{sec:3}

In this section we give inductive formulas for the kernels. Let us
start with the \pkernel. We set
$\p_2 : = \mu_2$ and
\begin{equation}
\label{Dny_otevrenych_dveri}
\p_n := \sum_{B} (-1)^{\vartheta(r_1,\ldots,r_k)}
\mu_k ( h \circ \p_{r_1} \otimes \cdots \otimes h \circ \p_{r_k}) 
\end{equation}
with the formal convention that $h \p_1 = \id$.  For our inductive
definition of the \qkernel\ we need the following notation:
\[
\p^i_n = \sum_{D} (-1)^{\vartheta(r_1,\ldots,r_{i-1})}
\mu_k(h \circ \p_{r_1} \otimes \cdots \otimes h \circ \p_{r_{i-1}}\ot
\otexp {\id}{n-i+1})
\]
where
\[
D := \{k,\Rada r1{i-1}\ |\
       2 \leq k \leq n,\ \ 
       i \leq k,\ r_1,\ldots, r_{i-1} \geq 1,\ 
       r_1+ \cdots + r_{i-1} + k - i +1 = n\},
\]
$i$ is a fixed integer, $1 \leq i \leq n$, and
where we again put $h\p_1 = \id_V$. 
We then define $\q_1 := \id_V$ and, inductively
\begin{equation}
\label{Eli_toto_utery_priletela}
\q_n = \sum_{C} (-1)^{n + r_i +\vartheta(r_1,\ldots,r_{i})}
\p^i_k (gf \circ \q_{r_1} \ot \cdots \ot gf \circ \q_{r_{i-1}} \ot h
\circ \q_{r_i} \ot \otexp {\id}{k-i})
\end{equation}

The first result of this note is:

\begin{theorem}
\label{main}
Let $\{\p_n\}_{n \geq 2}$ and $\{\q_n\}_{n \geq 1}$ be defined
inductively by~(\ref{Dny_otevrenych_dveri})
and~(\ref{Eli_toto_utery_priletela}). Then $\nu_n$, $\phi_n$, $\psi_n$
and $H_n$ determined by these $\p_n$ and $\q_n$ as in
formula~(\ref{eq:1}) solve Problem~\ref{problem}.
\end{theorem}

\noindent
{\bf Proof.}
A straightforward but awfully technical induction shows that the 
kernels satisfy:
\[
\delta (\p_n)
=
\sum_A
(-1)^{i(l+1) +n} \p_{k}(\id^{\ot i-1} \ot  gf \circ \p_l \ot
\id^{\ot k-i}),\ n \geq 2,
\]
and
\begin{eqnarray*}
\delta (\q_n) \hskip -2mm &=& \hskip -2mm -\sum_B
(-1)^{\vartheta(\Rada r1k)}
\p_k(\orada{gf \circ \q_{r_1}}{gf \circ \q_{r_k}})+
\\
\nonumber 
&&\hskip 1cm -
\sum_A
(-1)^{i(l+1) +n}
\q_{k}(\id^{\ot i-1} \ot \mu_l \ot \id^{\ot k-i}),\ n \geq 1.
\end{eqnarray*}
It is then almost obvious that the above two equations imply 
Axioms~(\ref{X})--(\ref{B}) for $\nu_n$, $\phi_n$, $\psi_n$ and $H_n$
defined by~(\ref{eq:1}).%
\qed

\section{Non-inductive formulas}
\label{sec:non-induct-form}
\label{sec:4}

In this section we give non-inductive formulas for the kernels. Our
formulas will be based on the language of trees which we use as names
for maps and their compositions. Formally this mean that we work in a
certain free operad, but we are not going to use this fancy language
here.  The terminology of trees is recalled in Section~II.1.5
of~\cite{markl-shnider-stasheff:book}.

\begin{figure}[t]
  \centering
\unitlength 1.3mm
\thicklines
\begin{picture}(60,36)(0,0)
\put(30,36){\line(0,-1){34}}
\put(30,32){\line(-1,-1){9.35}}
\put(19.45,21.45){\line(-1,-1){19.5}}
\put(20.1,21.9){\makebox(0,0)[cc]{\Large $\circ$}}
\put(9.8,12.2){\line(1,-1){2.8}}
\put(13.75,8.25){\line(1,-1){6.25}}
\put(13.25,8.7){\makebox(0,0)[cc]{\Large $\circ$}}
\put(30,32){\line(1,-1){9.35}}
\put(40.45,21.45){\line(1,-1){19.5}}
\put(39.9,21.9){\makebox(0,0)[cc]{\Large $\circ$}}
\put(16,6){\line(-1,-1){4}}
\put(50,12){\line(-1,-1){10}}
\put(50,12){\line(0,-1){10}}
\put(30,32){\makebox(0,0)[cc]{$\bullet$}}
\put(10,12){\makebox(0,0)[cc]{$\bullet$}}
\put(16,6){\makebox(0,0)[cc]{$\bullet$}}
\put(50,12){\makebox(0,0)[cc]{$\bullet$}}
\put(32,34){\makebox(0,0)[lb]{$u$}}
\put(9,14){\makebox(0,0)[rb]{$v$}}
\put(51,13){\makebox(0,0)[lb]{$w$}}
\put(19,7){\makebox(0,0)[rb]{$x$}}
\end{picture}
  \caption{An element of $\Trp_7$.}
  \label{example_of_tree}
\end{figure}
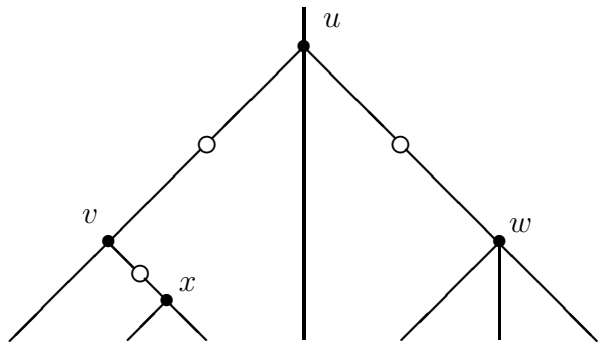

Let $\Trp_n$ denote the set of planar directed trees with at least binary
vertices (that is, all vertices have at least two incoming edges), 
with interior edges decorated by the symbol 
\krouzekstopka, and $n$ leaves. An example of such a tree is given
in Figure~\ref{example_of_tree}.
To each decorated tree $T \in \Trp_n$ we assign a map $F_T : \otexp Vn
\to V$, by interpreting $T$ as a ``flow chart,'' with
\krouzekstopka\  denoting the homotopy $h : V
\to V$ and a vertex of arity (= the number of incoming edges) $k$ denoting
the map $\mu_k : \otexp Vk \to V$.  For example, the tree $T$ in
Figure~\ref{example_of_tree} gives the degree $5$ map
\[
F_T = \mu_3(h \circ \mu_2(\id_V \ot h \circ \mu_2) \ot \id_V \ot h
\circ \mu_3) : \otexp V7 \to V
\]
which, evaluated at $(a,b,c,d,e,f,g) \in \otexp V7$, equals
\[
F_T(a,b,c,d,e,f,g) = 
(-1)^{|a|}\mu_3(h \circ 
\mu_2(a,h \circ \mu_2(b,c)),d,h\circ \mu_3(e,f,g)).  
\]

Finally, we assign to each tree $T \in \Trp_n$ the sign $\vartheta(T)$
as follows. For a vertex $v \in \Vert (T)$ of arity $k$ and $1 \leq i
\leq k$, let $r_i$ be the number of legs (= leaves) $e$ of $T$ such
that the unique path from $e$ to the root of $T$ contains the $i$-th
input edge of $v$. We then define $\vartheta_T(v) := \vartheta(\Rada
r1k)$ and $\vartheta(T) := \sum_{v \in \Vert(T)}\vartheta_T(v)$.

For example, for the tree $T$ in Figure~\ref{example_of_tree} we have,
at the vertex $u$ of arity $3$, $r_1=3$, $r_2 = 1$, $r_3 = 3$, at the
vertex $v$ of arity $2$, $r_1 =1$, $r_2 = 2$, at the vertex $w$ of
arity $3$, $r_1 = r_2 = r_3 = 1$ and, at the vertex $x$ of arity $2$,
$r_1 = r_2 = 1$. Therefore, modulo~$2$, $\vartheta_T(u) = 3\cdot 2 + 4
\cdot 4 = 0$, $\vartheta_T(v) = 1 \cdot 3 = 1$, $\vartheta_T(w) =
1\cdot 2 + 2 \cdot 2 = 0$ and $\vartheta_T(x) = 1 \cdot 2 = 0$, which
gives, again mod~$2$, $\vartheta(T) = 1$. We may finally formulate the
following almost obvious:

\begin{proposition}
\label{pred_odletem_do_Izraele}
The \pkernel\ $\p_n : \otexp Vn \to V$, defined inductively
by~(\ref{Dny_otevrenych_dveri}), can also be defined~as 
\[
\p_n := \sum_{T \in \Trp_n} (-1)^{\vartheta (T)} \cdot F_T,
\mbox { for each $n \geq 2$.}
\]
\end{proposition}

Let us proceed to our non-inductive definition of the \qkernel\ based
on a slightly more elaborate definition of decoration of a planar
tree.  We need to observe first that each planar tree $S$ admits a
natural total order of its set of vertices $\Vert(S)$ determined in the
following way.  

\begin{figure}[t]
\centering
\thicklines
{
\unitlength=.8pt
\begin{picture}(200.00,120.00)(0.00,0.00)
\put(170.00,50.00){\makebox(0.00,0.00){\scriptsize $1$}}
\put(90.00,30.00){\makebox(0.00,0.00){\scriptsize $2$}}
\put(10.00,30.00){\makebox(0.00,0.00){\scriptsize $3$}}
\put(40.00,60.00){\makebox(0.00,0.00){\scriptsize $4$}}
\put(110.00,110.00){\makebox(0.00,0.00){\scriptsize $5$}}
\put(80.00,20.00){\makebox(0.00,0.00){$\bullet$}}
\put(20.00,20.00){\makebox(0.00,0.00){$\bullet$}}
\put(160.00,40.00){\makebox(0.00,0.00){$\bullet$}}
\put(50.00,50.00){\makebox(0.00,0.00){$\bullet$}}
\put(100.00,100.00){\makebox(0.00,0.00){$\bullet$}}
\put(160.00,40.00){\line(0,-1){40.00}}
\put(160.00,40.00){\line(-1,-1){40.00}}
\put(100.00,100.00){\line(1,-1){100.00}}
\put(50.00,50.00){\line(1,-1){30.00}}
\put(80.00,20.00){\line(1,-1){20.00}}
\put(60.00,0.00){\line(1,1){20.00}}
\put(20.00,20.00){\line(1,-1){20.00}}
\put(100.00,100.00){\line(-1,-1){100.00}}
\put(100.00,120.00){\line(0,-1){20.00}}
\end{picture}}
\caption{\label{fig:3}Ordering vertices of a planar tree. The vertices
are numbered, from the biggest to the smallest one.}
\end{figure}
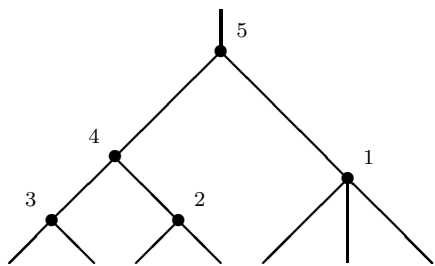

We say that a vertex $u$ is {\em below\/} a vertex
$v$ if $v$ lies on the (unique) directed path joining $u$ with the
root. This defines a {\em partial\/} order on the set of vertices of
$S$. It is easy to see that there exists precisely one {\em total\/} order $<$
on the set $\Vert(S)$ which satisfies the following two conditions:
\begin{itemize}
\item [(i)]
If $u$ is below $v$, then $u < v$.
\item [(ii)]
Suppose $S$ contains a subtree of the form shown in
Figure~\ref{pozitri_Izrael} and $1 \leq i \leq k-1$.
\begin{figure}[t]
\centering
\unitlength 2.3mm
\thicklines
\begin{picture}(31,17)(0,10)
\put(16,27){\line(0,-1){3}}
\put(16,24){\line(-6,-5){12}}
\put(4,14){\line(0,1){0}}
\put(16,24){\line(-1,-2){5}}
\put(16,24){\line(6,-5){12}}
\put(4,14){\line(-1,-1){3}}
\put(4,14){\line(-1,-3){1}}
\put(4,14){\line(1,-1){3}}
\put(11,14){\line(-1,-1){3}}
\put(11,14){\line(-1,-3){1}}
\put(11,14){\line(1,-1){3}}
\put(28,14){\line(-1,-1){3}}
\put(28,14){\line(-1,-3){1}}
\put(28,14){\line(1,-1){3}}
\put(16,24){\makebox(0,0)[cc]{$\bullet$}}
\put(4,14){\makebox(0,0)[cc]{$\bullet$}}
\put(11,14){\makebox(0,0)[cc]{$\bullet$}}
\put(28,14){\makebox(0,0)[cc]{$\bullet$}}
\put(4,15){\makebox(0,0)[rb]{$v_1$}}
\put(10,15){\makebox(0,0)[cb]{$v_2$}}
\put(28,15){\makebox(0,0)[lb]{$v_k$}}
\put(17,25){\makebox(0,0)[cc]{$v$}}
\put(17,17){\makebox(0,0)[cc]{$\cdots$}}
\put(5,11){\makebox(0,0)[cc]{$\cdots$}}
\put(12,11){\makebox(0,0)[cc]{$\cdots$}}
\put(29,11){\makebox(0,0)[cc]{$\cdots$}}
\end{picture}
\caption{\label{pozitri_Izrael}A subtree of $S$ used in the definition
  of the total order $<$.}
\end{figure}
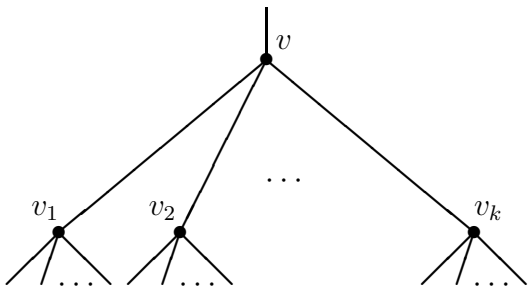
Then $v_i$ and all vertices below $v_i$ are less, in the order $<$,
than $v_{i+1}$.
\end{itemize}
See Figure~\ref{fig:3} for an example of such an order.

\begin{figure}[t]
\centering
{
\unitlength=1.000000pt
\begin{picture}(140.00,130.00)(0.00,0.00)
\thinlines
\put(120.00,10.00){\line(1,0){20.00}}
\put(60.00,30.00){\line(1,0){70.00}}
\put(10.00,50.00){\line(1,0){110.00}}
\put(30.00,80.00){\line(1,0){70.00}}
\put(50.00,100.00){\line(1,0){40.00}}
%
\thicklines
\put(130.00,20.00){\line(1,-2){10.00}}
\put(130.00,20.00){\line(-1,-2){10.00}}
\put(70.00,110.00){\line(2,-3){60.00}}
\put(70.00,40.00){\line(1,-2){20.00}}
\put(70.00,40.00){\line(-1,-2){20.00}}
\put(50.00,90.00){\line(2,-5){20.00}}
\put(20.00,60.00){\line(1,-3){20.00}}
\put(20.00,60.00){\line(-1,-3){20.00}}
\put(70.00,110.00){\line(-1,-1){50.00}}
\put(70.00,130.00){\line(0,-1){20.00}}
\put(130.00,20.00){\line(0,-1){20.00}}
\end{picture}}
\caption{\label{fig:4}%
Drawing horizontal lines.}
\end{figure}
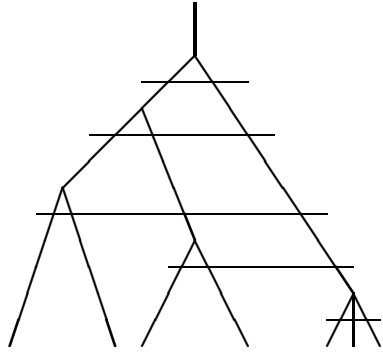

The next step is to redraw the tree in such a way that the vertices
are placed into different levels, according
to their order, and then draw horizontal lines slightly below the vertices,
as illustrated in Figure~\ref{fig:4}.
Now we decorate some (not necessary
all) of the intersections of the horizontal lines with the edges of
the tree with symbols \krouzekstopka\ or \teckastopka, according to
the following rules:
\begin{itemize}
\item[(i)] Let $\Rada x1k$ be the points at which a horizontal line
  intersects the edges of $S$, numbered from left to right. Then there
  is some $0 \leq s \leq k-1$ such that the points $\Rada x1s$ are
  decorated by \teckastopka, $x_{s+1}$ is decorated by \krouzekstopka\
  and the points $\Rada x{s+2}k$ are not decorated.
\item[(ii)]
Each edge of $S$ is decorated at most once.
\item[(iii)]
Each internal edge of $S$ is decorated.
\end{itemize}
 
Condition~(i) means that we may see the
following pattern%
\footnote{This should remind us about the time when this
  paper was finished -- carp with potato salad is the most typical
  Czech Christmas dish.}
on the horizontal lines:
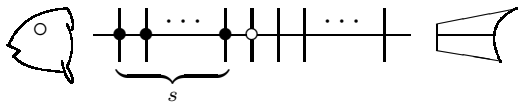
\begin{center}
{
\unitlength=1.000000pt
\label{fish}
\begin{picture}(120.00,25.00)(0.00,0.00)
\put(-40,-10){\thinlines\unitlength=2.000000pt
\qbezier(4,8)(6.5,14.5)(15,15)
\qbezier(6,7)(4,5)(6,7)
\qbezier(6,7)(4.5,5.5)(5,6)
\qbezier(5,6)(6.5,3.5)(14,3)
\qbezier(15,15)(16.5,13)(16,11)
\qbezier(15,12)(17,8)(15,4)
\qbezier(15,4)(16.5,1.5)(16,1)
\qbezier(16,1)(15,1)(14,5)
\put(10,11){\makebox(0,0)[cc]{$\circ$}}
\qbezier(4,8)(3,8.5)(6,7)
}
\put(100,-34){\thinlines\unitlength=2pt
\put(15,24){\line(0,-1){5}}
\put(15,24){\line(5,1){15}}
\put(15,19){\line(6,-1){12}}
\qbezier(30,27)(23.5,24)(27,17)
\put(27,25){\line(0,1){0}}
\put(15,22){\line(1,0){11}}
}
\thinlines
\put(94.00,15.00){\makebox(0.00,0.00){$\cdots$}}
\put(34.00,15.00){\makebox(0.00,0.00){$\cdots$}}
\put(60.00,10.00){\makebox(0.00,0.00){$\circ$}}
\put(50.00,10.00){\makebox(0.00,0.00){$\bullet$}}
\put(20.00,10.00){\makebox(0.00,0.00){$\bullet$}}
\put(10.00,10.00){\makebox(0.00,0.00){$\bullet$}}
\put(0.00,10.00){\line(1,0){58}}
\put(62,10.00){\line(1,0){58}}
\thicklines
\put(110.00,20.00){\line(0,-1){20.00}}
\put(80.00,20.00){\line(0,-1){20.00}}
\put(70.00,20.00){\line(0,-1){20.00}}
\put(60.00,0.00){\line(0,1){7.9}}
\put(60.00,20.00){\line(0,-1){7.5}}
\put(50.00,20.00){\line(0,-1){20.00}}
\put(20.00,20.00){\line(0,-1){20.00}}
\put(10.00,20.00){\line(0,-1){20.00}}
\put(30,-2){\makebox(0.00,0.00)[t]{$\underbrace{\rule{15mm}{0mm}}_s$}}
\end{picture}}
\end{center}

\vskip 3mm
\noindent 
with the case $s=0$ (no black dot) allowed. 

\begin{figure}[t]
\centering
{
\unitlength=1.000000pt
\begin{picture}(140.00,130.00)(0.00,0.00)
\thinlines
\qbezier(30.00,80.00)(33.75,79.75)(37.75,80.00)
\put(100.00,80.00){\line(-1,0){58.00}}
\qbezier(120.00,10.00)(122.00,10)(123.25,10)
\put(140.00,10.00){\line(-1,0){12.5}}
\put(50.00,100.00){\line(1,0){25.0}}
\put(90.00,100.00){\line(-1,0){11}}
\put(120.00,50.00){\line(-1,0){52.0}}
\put(10.25,50.00){\line(1,0){53.50}}
\put(130.00,30.00){\line(-1,0){53.00}}
\put(73.00,30.00){\line(-1,0){13.00}}
\thicklines
\put(20.00,60.00){\line(1,1){18.25}}
\put(70.00,110.00){\line(-1,-1){28.75}}
\qbezier(129.00,20.00)(127.50,16.25)(126.0,12.0)
\qbezier(121.00,0.25)(123,5.25)(124,8.0)
\qbezier(70.00,110.00)(72.50,106.25)(75.00,102.00)
\put(130.00,20.00){\line(-2,3){52}}
\put(50,0.00){\line(1,2){20.00}}
\put(90.00,0.00){\line(-1,2){14}}
\put(74.50,32.50){\line(-1,2){7.75}}
\put(50.00,90.00){\line(2,-5){15.00}}
\put(130.00,20.00){\line(1,-2){10.00}}
\put(20.00,60.00){\line(1,-3){20.00}}
\put(20.00,60.00){\line(-1,-3){20.00}}
\put(70.00,130.00){\line(0,-1){20.00}}
\put(125.50,9.75){\makebox(0.00,0.00){$\circ$}}
\put(65.25,30.00){\makebox(0.00,0.00){$\bullet$}}
\put(66,50.00){\makebox(0.00,0.00){$\circ$}}
\put(77.25,99.75){\makebox(0.00,0.00){$\circ$}}
\put(40.25,79.75){\makebox(0.00,0.00){$\circ$}}
\put(75.00,30){\makebox(0.00,0.00){$\circ$}}
\put(23.75,50.00){\makebox(0.00,0.00){$\bullet$}}
\put(16.75,50.00){\makebox(0.00,0.00){$\bullet$}}
\put(60.00,100.00){\makebox(0.00,0.00){$\bullet$}}
\put(130.00,20.00){\line(0,-1){20.00}}
%
\end{picture}}
\caption{\label{fig:4a}%
A decoration of the tree from the previous figure.}
\end{figure}
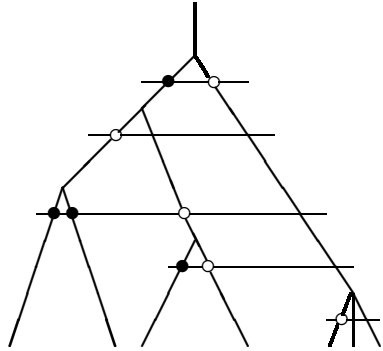

A decoration of the tree
from Figure~\ref{fig:4} is shown in Figure~\ref{fig:4a}. 
All possible decorations of the tree \hskip .5mm
{
\unitlength=.15pt
\begin{picture}(100.00,70.00)(0.00,0.00)
\qbezier(20.00,20.00)(30.00,10.00)(40.00,0.00)
\qbezier(80.00,20.00)(70.00,10.00)(60.00,0.00)
\qbezier(50.00,50.00)(80.00,20.00)(100.00,0.00)
\qbezier(50.00,50.00)(20.00,20.00)(0.00,0.00)
\qbezier(50.00,70.00)(50.00,60.00)(50.00,50.00)
\end{picture}}
\hskip 1.5mm are shown in Figure~\ref{fig:5}.

\begin{figure}[t]
  \centering
\unitlength .59mm
\ifx\plotpoint\undefined\newsavebox{\plotpoint}\fi 
\begin{picture}(72,95)(0,-45)
\thicklines
\put(-80,0){
\put(8,7){\line(1,1){6}}
\put(14,13){\line(0,1){0}}
\put(16,15){\line(1,1){14}}
\put(32,31){\line(1,1){8}}
\put(40,39){\line(0,1){0}}
\put(40,39){\line(1,-1){8}}
\put(48,31){\line(0,1){0}}
\put(50,29){\line(1,-1){14}}
\put(66,13){\line(1,-1){6}}
\put(22,21){\line(1,-1){6}}
\put(30,13){\line(1,-1){6}}
\put(58,21){\line(-1,-1){6}}
\put(50,13){\line(-1,-1){6}}
\put(40,45){\line(0,-1){6}}
\put(31,30){\makebox(0,0)[cc]{$\circ$}}
\put(49,30){\makebox(0,0)[cc]{$\circ$}}
\put(15,14){\makebox(0,0)[cc]{$\bullet$}}
\put(29,14){\makebox(0,0)[cc]{$\bullet$}}
\put(51,14){\makebox(0,0)[cc]{$\circ$}}
\put(65,14){\makebox(0,0)[cc]{\line(1,-1){10}}}
}
\put(0,0){
\put(8,7){\line(1,1){6}}
\put(14,13){\line(0,1){0}}
\put(16,15){\line(1,1){14}}
\put(32,31){\line(1,1){8}}
\put(40,39){\line(0,1){0}}
\put(40,39){\line(1,-1){8}}
\put(48,31){\line(0,1){0}}
\put(50,29){\line(1,-1){14}}
\put(66,13){\line(1,-1){6}}
\put(22,21){\line(1,-1){6}}
\put(30,13){\line(1,-1){6}}
\put(58,21){\line(-1,-1){6}}
\put(50,13){\line(-1,-1){6}}
\put(40,45){\line(0,-1){6}}
\put(31,30){\makebox(0,0)[cc]{$\circ$}}
\put(49,30){\makebox(0,0)[cc]{$\circ$}}
\put(15,14){\makebox(0,0)[cc]{$\bullet$}}
\put(29,14){\makebox(0,0)[cc]{$\bullet$}}
\put(51,14){\makebox(0,0)[cc]{$\bullet$}}
\put(65,14){\makebox(0,0)[cc]{$\circ$}}
}
%
\put(80,0){
\put(8,7){\line(1,1){6}}
\put(14,13){\line(0,1){0}}
\put(16,15){\line(1,1){14}}
\put(32,31){\line(1,1){8}}
\put(40,39){\line(0,1){0}}
\put(40,39){\line(1,-1){8}}
\put(48,31){\line(0,1){0}}
\put(50,29){\line(1,-1){14}}
\put(66,13){\line(1,-1){6}}
\put(22,21){\line(1,-1){6}}
\put(30,13){\line(1,-1){6}}
\put(58,21){\line(-1,-1){6}}
\put(50,13){\line(-1,-1){6}}
\put(40,45){\line(0,-1){6}}
\put(31,30){\makebox(0,0)[cc]{$\bullet$}}
\put(49,30){\makebox(0,0)[cc]{$\circ$}}
\put(15,14){\makebox(0,0)[cc]{$\circ$}}
\put(29,14){\makebox(0,0)[cc]{\line(1,-1){10}}}
\put(51,14){\makebox(0,0)[cc]{$\circ$}}
\put(65,14){\makebox(0,0)[cc]{\line(1,-1){10}}}
}
\put(0,-50){
\put(-80,0){
\put(8,7){\line(1,1){6}}
\put(14,13){\line(0,1){0}}
\put(16,15){\line(1,1){14}}
\put(32,31){\line(1,1){8}}
\put(40,39){\line(0,1){0}}
\put(40,39){\line(1,-1){8}}
\put(48,31){\line(0,1){0}}
\put(50,29){\line(1,-1){14}}
\put(66,13){\line(1,-1){6}}
\put(22,21){\line(1,-1){6}}
\put(30,13){\line(1,-1){6}}
\put(58,21){\line(-1,-1){6}}
\put(50,13){\line(-1,-1){6}}
\put(40,45){\line(0,-1){6}}
\put(31,30){\makebox(0,0)[cc]{$\bullet$}}
\put(49,30){\makebox(0,0)[cc]{$\circ$}}
\put(15,14){\makebox(0,0)[cc]{$\bullet$}}
\put(29,14){\makebox(0,0)[cc]{$\circ$}}
\put(51,14){\makebox(0,0)[cc]{$\circ$}}
\put(65,14){\makebox(0,0)[cc]{\line(1,-1){10}}}
}
\put(0,0){
\put(8,7){\line(1,1){6}}
\put(14,13){\line(0,1){0}}
\put(16,15){\line(1,1){14}}
\put(32,31){\line(1,1){8}}
\put(40,39){\line(0,1){0}}
\put(40,39){\line(1,-1){8}}
\put(48,31){\line(0,1){0}}
\put(50,29){\line(1,-1){14}}
\put(66,13){\line(1,-1){6}}
\put(22,21){\line(1,-1){6}}
\put(30,13){\line(1,-1){6}}
\put(58,21){\line(-1,-1){6}}
\put(50,13){\line(-1,-1){6}}
\put(40,45){\line(0,-1){6}}
\put(31,30){\makebox(0,0)[cc]{$\bullet$}}
\put(49,30){\makebox(0,0)[cc]{$\circ$}}
\put(15,14){\makebox(0,0)[cc]{$\circ$}}
\put(29,14){\makebox(0,0)[cc]{\line(1,-1){10}}}
\put(51,14){\makebox(0,0)[cc]{$\bullet$}}
\put(65,14){\makebox(0,0)[cc]{$\circ$}}
}
%
\put(80,0){
\put(8,7){\line(1,1){6}}
\put(14,13){\line(0,1){0}}
\put(16,15){\line(1,1){14}}
\put(32,31){\line(1,1){8}}
\put(40,39){\line(0,1){0}}
\put(40,39){\line(1,-1){8}}
\put(48,31){\line(0,1){0}}
\put(50,29){\line(1,-1){14}}
\put(66,13){\line(1,-1){6}}
\put(22,21){\line(1,-1){6}}
\put(30,13){\line(1,-1){6}}
\put(58,21){\line(-1,-1){6}}
\put(50,13){\line(-1,-1){6}}
\put(40,45){\line(0,-1){6}}
\put(31,30){\makebox(0,0)[cc]{$\bullet$}}
\put(49,30){\makebox(0,0)[cc]{$\circ$}}
\put(15,14){\makebox(0,0)[cc]{$\bullet$}}
\put(29,14){\makebox(0,0)[cc]{$\circ$}}
\put(51,14){\makebox(0,0)[cc]{$\bullet$}}
\put(65,14){\makebox(0,0)[cc]{$\circ$}}
}}
\end{picture}
   \caption{All possible decorations of a tree.\label{fig:5}}
\end{figure}
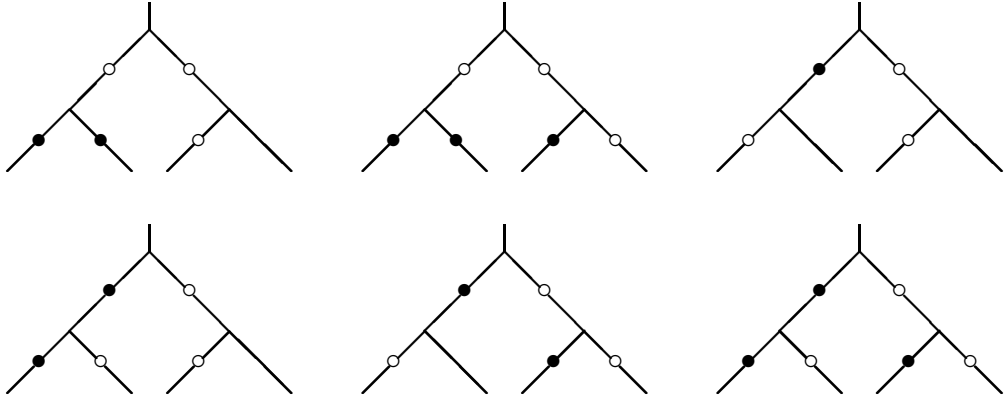

Let $\Trq_n$ be the set of all decorated, in the above sense, planar
directed trees with at least binary vertices and $n$ leaves. To each
$S \in \Trq_n$ we assign a map $G_S : \otexp Vn \to V$, by interpreting
$S$ as a ``flow chart,'' with \krouzekstopka\ denoting the homotopy $h
: V \to V$, \teckastopka\ denoting the composition $gf$, and a vertex
of arity $k$ the map $\mu_k : \otexp Vk \to V$.  For example, the tree
$S$ in Figure~\ref{fig:4a} gives degree $6$ map
\[
\mu_2(gf \circ \mu_2(h \circ \mu_2(gf \ot gf)\ot h \circ \mu_2 (gf \ot
h))\ot h \circ \mu_3 (h \ot \otexp {\id}2)) : \otexp V7 \to V.
\]

Finally, we must define a sign $\varepsilon(S)$ of a tree $S \in
\Trq_n$. The definition is more difficult than the definition of the sign
$\vartheta(T)$ of a tree $T \in \Trp_n$, because 
$\varepsilon(S)$ will depend also on the decoration, not only on the
combinatorial type, of the tree $S$. 

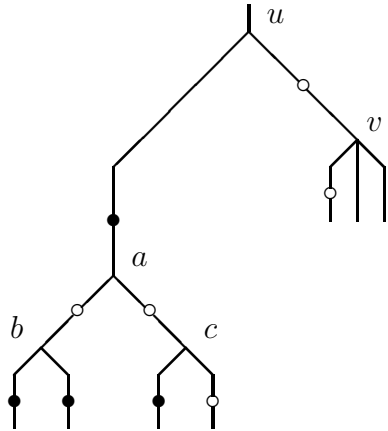
\begin{figure}[t]
  \centering
\unitlength 1.2mm
\thicklines
\begin{picture}(64,50)(10,0)
\put(46,47){\line(0,-1){3}}
\put(46,44){\line(1,-1){5.5}}
\put(52.4,37.6){\line(1,-1){5.5}}
\put(58,32){\line(-1,-1){3}}
\put(58,32){\line(0,-1){3}}
\put(58,32){\line(1,-1){3}}
\put(46,44){\line(-1,-1){15}}
\put(31,29){\line(0,-1){6}}
\put(55,29){\line(0,-1){2.4}}
\put(55,24){\line(0,-1){1}}
\put(55,25.5){\line(0,-1){2}}
\put(58,29){\line(0,-1){6}}
\put(61,29){\line(0,-1){6}}
\put(31,23){\line(0,1){0}}
\put(31,23){\line(0,-1){6}}
\put(31,17){\line(0,1){0}}
\put(31,17){\line(-1,-1){3.5}}
\put(26.5,12.5){\line(-1,-1){3.5}}
\put(31,17){\line(1,-1){3.5}}
\put(35.5,12.5){\line(1,-1){3.5}}
\put(23,9){\line(-1,-1){3}}
\put(23,9){\line(1,-1){3}}
\put(39,9){\line(-1,-1){3}}
\put(39,9){\line(1,-1){3}}
\put(20,6){\line(0,-1){6}}
\put(26,6){\line(0,-1){6}}
\put(36,6){\line(0,-1){2}}
\put(36,2){\line(0,-1){2}}
\put(36,5){\line(0,-1){3}}
\put(42,6){\line(0,-1){2.4}}
\put(42,1){\line(0,1){1.4}}
\put(42,1){\line(0,-1){1}}
\put(48,45){\makebox(0,0)[bl]{$u$}}
\put(59,33){\makebox(0,0)[bl]{$v$}}
\put(31,23){\makebox(0,0)[cc]{$\bullet$}}
\put(33,18){\makebox(0,0)[bl]{$a$}}
\put(21,10){\makebox(0,0)[br]{$b$}}
\put(41,10){\makebox(0,0)[bl]{$c$}}
\put(20,3){\makebox(0,0)[cc]{$\bullet$}}
\put(26,3){\makebox(0,0)[cc]{$\bullet$}}
\put(36,3){\makebox(0,0)[cc]{$\bullet$}}
\put(42,3){\makebox(0,0)[cc]{$\circ$}}
\put(27,13){\makebox(0,0)[cc]{$\circ$}}
\put(35,13){\makebox(0,0)[cc]{$\circ$}}
\put(55,26){\makebox(0,0)[cc]{$\circ$}}
\put(52,38){\makebox(0,0)[cc]{$\circ$}}
\end{picture}
\caption{Decomposing a tree into \pkernels.}
\label{fig:6}
\end{figure}

To calculate $\varepsilon(S)$, we must first decompose $S$ into trees
$T_1,\ldots,T_k$ representing summands of \pkernels, following the
pattern of~(\ref{Eli_toto_utery_priletela}). The sign is then defined
as 
\[
\varepsilon(S) :=  n + r_i +\vartheta(r_1,\ldots,r_i)
+ \sum_1^k \vartheta(T_j),
\]
where $\Rada r1i$ have the same meaning as
in~(\ref{Eli_toto_utery_priletela}). Let us calculate, as an example,
the sign of the decorated tree from Figure~\ref{fig:4a}. The
decomposition of this tree into trees from $\Trp$ is shown in
Figure~\ref{fig:6}.
In this figure, $T_1$ is the decorated subtree with vertices $u$ and
$v$ and $T_2$ is the subtree with vertices $a, b$ and $c$. The sign of
$S$ is then the sum $\varepsilon(S) := 7+3+ 1\cdot 4 + \vartheta(T_1)
+ \vartheta(T_2) = 0 \mbox{ (mod~2)}$.  The following proposition then
follows from a boring combinatorics.

\begin{proposition}
\label{po_Xmas_party_na_MFF_mne_boli_siska}
The \qkernel\ $\q_n : \otexp Vn \to V$, defined inductively
by~(\ref{Eli_toto_utery_priletela}), can also be defined~as 
\[
\q_n := \sum_{S \in \Trq_n} (-1)^{\varepsilon (S)} \cdot G_S,
\mbox { for $n \geq 2$,}
\]
and $q_1 := \id_V$ for $n=1$.
\end{proposition}

\section{Why do the transfers exist}
\label{sec:5}

As we already observed, if the basic ring $R$ is a field of in
characteristic $0$, the existence of transfers follows from a general
theory developed in~\cite{markl:ha} -- see ``move {\bf (S)}'' on page
141 of~\cite{markl:ha}. We want to make this statement more precise now.
In this section we assume that the reader is familiar with colored
operads which describe diagrams of algebras, 
see~\cite{markl:ha} again. The rest of the paper is independent on the
material of this section.

Let $\Pin$ be the $2$-colored operad describing structures consisting
of an associative multiplication $\mu$ on a vector space 
$V$ and linear maps of vector spaces $f : V
\to W$, $g : W \to V$ such that $gf = \id_V$.  Let also $\Pout$ be the
$2$-colored operad describing diagrams consisting of an associative
multiplication $\mu$ on $V$, an associative multiplication $\nu$ on
$W$, and homomorphisms $f : V \to W$, $g : W \to V$ of these
associative algebras such that $gf = \id_V$. An explicit description
of these operads can be found in Example~12 of~\cite{markl:ha}, where
$\Pin$ was denoted ${\cal P}_{(S,\underline D)}$ and $\Pout$ was
denoted ${\cal P}_{\underline S}$.  Finally, let $S : \Pout \to \Pin$
be the map defined by
\[
S(\mu) := \mu,\ S(f) := f,\ S(g) := g\ \mbox { and }\
S(\nu) := f\mu ( g \ot g).   
\]
This well-defined map of colored operads represents a solution of the
following ``classical limit'' of Problem~\ref{problem}.

\begin{problem}
\label{class}
We are given two vector spaces $V$, $W$ and linear maps $f:V \to W$,
$g: W \to V$ such that $gf = \id_V$ (in other words, $f: V
\hookrightarrow W$ is an inclusion and $g$ its retraction). Given an
associative algebra structure $\mu : V \ot V \to V$ on the vector
space $V$, find an associative algebra structure $\nu : W \ot W \to
W$ on $W$ such that $f$ and $g$ became homomorphisms of associative
algebras.
\end{problem}

Let $\Rin$ be the dg-operad representing the ``input data'' of our
transfer problem for $\ainfty$-algebras, that is, diagrams consisting
of an $\ainfty$-structure $\bfmu = (\mu_2,\mu_3,\ldots)$ on
$(V,\pa_V)$, dg-maps $f : (V,\pa_V) \to (W,\pa_W)$, $g : (W,\pa_W) \to
(V,\pa_V)$ and a chain homotopy $h$ between $gf$ and $\id_V$.  Let
$\rhoin : \Rin \to \Pin$ be the map of colored operad given by
\[
\rhoin(\mu_2) := \mu,\ \rhoin(\mu_n) := 0\ \mbox { for $n \geq 3$},\
\rhoin(f) := f,\ \rhoin(g) := g\ \mbox { and }\ \rhoin(h) := 0.
\]

Similarly, let $\Rout$ be the dg-operad representing a solution of
our transfer problem, that is, diagrams consisting of an
$\ainfty$-structure $\bfmu = (\mu_2,\mu_3,\ldots)$ on $(V,\pa_V)$, an
$\ainfty$-structure $\bfnu = (\nu_2,\nu_3,\ldots)$ on $(W,\pa_W)$,
$\ainfty$-maps $\bfphi = (\phi_1,\phi_2,\ldots) : (V,\pa,\bfmu) \to
(W,\pa,\bfnu)$, $\bfpsi = (\psi_1,\psi_2,\ldots) : (W,\pa,\bfnu)
\to (V,\pa,\bfmu)$ and an $\ainfty$-homotopy $\bfH = (H_1,H_2,\ldots)$
between $\bfpsi\bfphi$ and $\id_V$.  Let $\rhoout : \Rout \to \Pout$
be the map defined by
\begin{eqnarray*}
&
\rhoout(\mu_2) := \mu,\ \rhoout(\mu_n) := 0\ \mbox { for $n \geq 3$, }\
\rhoout(\nu_2) := \nu,\ \rhoout(\nu_n) := 0\ \mbox { for $n \geq 3$,}
&\\&
\rhoout(\phi_1) := f,\ \rhoout(\phi_n) := 0\ \mbox { for $n \geq 2$, }\
\rhoout(\psi_1) := g,\ \rhoout(\psi_n) := 0\ \mbox { for $n \geq 2$, and}\
&\\&
\rhoout(H_n) := 0\ \mbox { for $n \geq 1$.}
&
\end{eqnarray*}
The following proposition follows from the methods 
of~\cite{markl:ha} and~\cite{markl:ho}.

\begin{proposition}
The map $\rhoin : \Rin \to \Pin$ is a cofibrant resolution of the
colored operad $\Pin$ and $\rhoout : \Rout \to \Pout$ is a cofibrant
resolution of the colored operad $\Pout$.
\end{proposition}

It follows from~\cite[Lemma~20]{markl:ha} 
that there exist a lift $\widetilde S :
\Rout \to \Rin$ making the following diagram commutative:
\[
\square {\Rout}{\Rin}{\Pout}{\Pin}{\widetilde S}{\rhoin}{\rhoout}{S}
\]
Each such a lift $\widetilde S$ clearly provides a solution of
Problem~\ref{problem} while formulas~(\ref{eq:1}) determine a specific
lift of $S$.

Observe that, very crucially, we work with algebras {\em without
units\/}. It is straightforward realize that the unital version
of Problem~\ref{class} does not have an affirmative answer.  Indeed,
given a unit $1_V \in V$ for $\mu$, we would be forced to define $1_W
:= f(1_V)$. It is then easy to see that such $1_W$ is a unit for the
transferred structure if and only if $fg = \id_W$, that is, $f$ and
$g$ are isomorphisms inverse to each other which we did not assume.
This observation explains why our solution of Problem~\ref{problem}
which is, as we explained above, a lift of the classical
Problem~\ref{class}, works only for non-unital $A_\infty$-algebras.
Transfers of unital $A_\infty$-structures present much harder problem,
see the analysis in~\cite{markl:ha}.

\section{Some other properties of the transfer}
\label{sec:6}

In this section we analyze what happens if $g$ is not just a left
homotopy inverse of $f$, but if $f$ and $g$ are chain homotopy
equivalences inverse to each other.

Let $\setinf(V,\pa)$ denote the set of {\em isomorphism classes\/}
(with respect to $\ainfty$-maps) of $\ainfty$-structures on a given
chain complex $(V,\pa)$. Suppose we are given chain maps $f :
(V,\pa_V) \to (W,\pa_W)$, $g : (W,\pa_W) \to (V,\pa_V)$ and a chain
homotopy $h$ between $gf$ and $\id_V$. It can be easily shown that the first
formula of~(\ref{eq:1}) defines a set map
\[
\Trans fgh : \setinf(V,\pa_V) \to \setinf(W,\pa_W).
\]
Suppose we are given also a chain homotopy $l$ between $fg$ and
$\id_W$, that is, $f$ and $g$ are now fully fledged chain homotopy
equivalences inverse to each other. Then one may as well consider the
map
\[
\Trans gfl : \setinf(W,\pa_W) \to \setinf(V,\pa_V).
\]
We found the following proposition surprising, because there is no
relation between the homotopies $h$ and $l$.

\begin{proposition}
\label{co_zitra?}
Let $f$ and $g$ be chain homotopy equivalences, with chain homotopies
$h : gf \cong \id_V$ and $l : fg \cong \id_W$. Then both $\Trans fgh$
and $\Trans gfl$ are isomorphisms and
\[
\Trans gfl = \Trans fgh^{-1}.
\]
\end{proposition}

\noindent 
{\bf Proof.}  Formulas~(\ref{eq:1}) give an $\ainfty$-structure
$\bfnu$ on $(W,\pa_W)$ and an $\ainfty$-map $\bfphi : (V,\pa_V,\bfnu)
\to (W,\pa_W,\bfmu)$. Let us apply~(\ref{eq:1}) once again, this time
to construct an $\ainfty$-structure $\bar{\bfmu}$ on $(V,\pa_V)$
together with an $\ainfty$-map $\bar{\bfphi} : (W,\pa_W,\bfnu) \to
(V,\pa_V,\bar{\bfmu})$, using $g$ instead of $f$, $f$
instead of $g$ and $l$ instead of $h$. We must prove that $\bfmu$ is
isomorphic to $\bar{\bfmu}$. To this end, recall the following
$\ainfty$-case of ``move {\bf (M2)}'' of~\cite{markl:ha}.

\begin{proposition}
\label{uz_si_na_ni_nekdy_ani_nevzpomenu}
Let $(A,\pa_A,\bfxi)$, $(B,\pa_B,\bfeta)$ be $\ainfty$-algebras and
$\bftheta = (\theta_1,\theta_2,\ldots) : (A,\pa_A,\bfxi)\to 
(B,\pa_B,\bfeta)$ an $\ainfty$-map. Suppose that
$C : (A,\pa_A) \to (B,\pa_B)$ is a chain map, homotopic to the linear
part $\theta_1$ of $\bftheta$. Then $C$ can be extended into an
$\ainfty$-map $\bfC = (C_1 = C, C_2,\ldots) :  (A,\pa_A,\bfxi)\to 
(B,\pa_B,\bfeta)$. 
\end{proposition}

Now observe that the linear part of the composition
$\bar{\bfphi}\bfphi$ equals
$\id_V$. Proposition~\ref{uz_si_na_ni_nekdy_ani_nevzpomenu} then implies
the existence of an $\ainfty$-map $\bfC = (\id_V,C_2,\ldots):
(V,\pa_V,\bfmu) \to (V,\pa_V,\bar{\bfmu})$ which is clearly an
isomorphism. This finishes our proof of Proposition~\ref{co_zitra?}.
Observe that the composition $\bar{\bfphi}\bfphi$ need not be an
isomorphism, therefore the full force of
Proposition~\ref{uz_si_na_ni_nekdy_ani_nevzpomenu} is necessary.%
\qed

Let us consider again chain homotopy equivalences $f$ and $g$, with
chain homotopies $h : gf \cong \id_V$ and $l : fg \cong \id_W$.  Given
an $\ainfty$-structure $\bfmu = (\mu_2,\mu_3,\ldots)$ on $(V,\pa_V)$,
let us construct, using formulas~(\ref{eq:1}), an $\ainfty$-structure
$\bfnu = (\nu_1,\nu_2,\ldots)$ on $(W,\pa_W)$ and $\ainfty$-maps
$\bfphi$, $\bfpsi$ as before.
A natural question is when such a situation gives rise to a
``perfect'' chain homotopy equivalence in the category of
$\ainfty$-algebras. The following proposition follows from the methods
of~\cite{markl:ip}. 

\begin{proposition}
The chain homotopy $l$ can be extended into an $\ainfty$-homotopy
$\bfL$ between $\ainfty$-maps $\bfphi \bfpsi$ if the chain
homotopy equivalence $(f,g,h,l)$ extends into a strong homotopy
equivalence in the sense of~\cite[Definition~1]{markl:ip}. This,
according to~\cite[Theorem~11]{markl:ip}, happens if and only if
\begin{equation}
\label{posledni_rovnice}
\hskip 3mm
[fh-lf] = 0 \mbox { in } H_1(\Hom(V,W)) \mbox { or, equivalently, }
[gl-hg] = 0 \mbox { in } H_1(\Hom(W,V)).
\end{equation}
\end{proposition}

If the $\ainfty$ structure $\bfmu = (\mu_2,\mu_3,\ldots)$ on
$(V,\pa_V)$ is generic enough, then the vanishing of the obstruction
classes in~(\ref{posledni_rovnice}) is also necessary for the
existence of an extension of $l$ into $\bfL$.

\section{Two observations}
\label{sec:7}

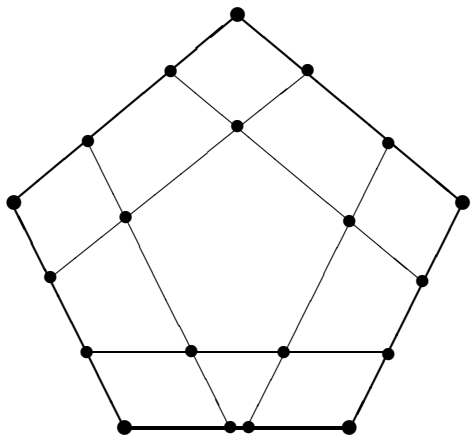
\begin{figure}[t]
  \centering
\unitlength 5mm
\begin{picture}(10,11)(10,8)
\thicklines
\put(8,14){\line(6,5){6}}
\put(14,19){\line(6,-5){6}}
\put(16,8){\line(-1,0){5}}
\put(11,8){\line(-1,2){3}}
\put(20,14){\line(-1,-2){3}}
\put(17,8){\line(-1,0){1}}
\thinlines
\put(8.993,11.966){\line(5,4){6.875}}
\put(9.997,10.034){\line(1,0){7.916}}
\put(18.061,15.608){\line(-1,-2){3.846}}
\put(12.189,17.504){\line(6,-5){6.689}}
\put(9.997,15.645){\line(1,-2){3.791}}
\put(14.01,18.99){\makebox(0,0)[cc]{\large$ \bullet$}}
\put(12.227,17.466){\makebox(0,0)[cc]{$\bullet$}}
\put(10.034,15.608){\makebox(0,0)[cc]{$\bullet$}}
\put(8.064,13.973){\makebox(0,0)[cc]{\large $\bullet$}}
\put(9.031,12.004){\makebox(0,0)[cc]{$\bullet$}}
\put(11.037,13.602){\makebox(0,0)[cc]{$\bullet$}}
\put(14.01,16.017){\makebox(0,0)[cc]{$\bullet$}}
\put(15.868,17.504){\makebox(0,0)[cc]{$\bullet$}}
\put(18.024,15.571){\makebox(0,0)[cc]{$\bullet$}}
\put(16.983,13.49){\makebox(0,0)[cc]{$\bullet$}}
\put(19.994,13.973){\makebox(0,0)[cc]{\large $\bullet$}}
\put(18.916,11.892){\makebox(0,0)[cc]{$\bullet$}}
\put(9.997,9.997){\makebox(0,0)[cc]{$\bullet$}}
\put(12.784,10.034){\makebox(0,0)[cc]{$\bullet$}}
\put(15.237,9.997){\makebox(0,0)[cc]{$\bullet$}}
\put(18.024,9.96){\makebox(0,0)[cc]{$\bullet$}}
\put(11,8){\makebox(0,0)[cc]{\large $\bullet$}}
\put(13.825,8){\makebox(0,0)[cc]{$\bullet$}}
\put(14.308,8){\makebox(0,0)[cc]{$\bullet$}}
\put(16.983,8){\makebox(0,0)[cc]{\large $\bullet$}}
\end{picture}
 \caption{The decomposition of the associahedron $K_4$ induced by
  $\p_4$. It consists of 10 squares and one pentagon. The squares
  adjacent to the vertices of $K_4$ correspond to the five trees of $\Trp_4$
  with two interior edges, the squares adjacent to the edges of $K_4$
  correspond to the five trees of $\Trp_4$ with one interior edge.
  The pentagon in the center of $K_4$ corresponds to the corrola (tree
  with no interior edge) in $\Trp_4$.}
  \label{fig:ob3}
\end{figure}

\begin{figure}[t]
  \centering
\unitlength 3mm
\begin{picture}(24.2,17.021)(1,2)
%
%
\put(-15,0){
\thicklines
\put(10,17){\line(1,0){10}}
\put(20,17){\line(3,-5){4.2}}
\put(24.2,10){\line(-3,-5){4.2}}
\put(20,3){\line(-1,0){10}}
\put(10,3){\line(-3,5){4.2}}
\put(5.8,10){\line(3,5){4.2}}
\thinlines
\put(24,10){\line(-3,1){15}}
\put(20,3){\line(-3,2){13}}
\put(10.108,17.021){\makebox(0,0)[cc]{\large $\bullet$}}
\put(19.994,17.021){\makebox(0,0)[cc]{\large $\bullet$}}
\put(5.797,9.997){\makebox(0,0)[cc]{\large $\bullet$}}
\put(9.997,2.973){\makebox(0,0)[cc]{\large $\bullet$}}
\put(19.994,2.973){\makebox(0,0)[cc]{\large $\bullet$}}
\put(24.156,9.96){\makebox(0,0)[cc]{\large $\bullet$}}
\put(8.882,15.051){\makebox(0,0)[cc]{$\bullet$}}
\put(6.875,11.706){\makebox(0,0)[cc]{$\bullet$}}
}
%
\put(10,0){
\thicklines
\put(10,17){\line(1,0){10}}
\put(20,17){\line(3,-5){4.2}}
\put(24.2,10){\line(-3,-5){4.2}}
\put(20,3){\line(-1,0){10}}
\put(10,3){\line(-3,5){4.2}}
\put(5.8,10){\line(3,5){4.2}}
\put(10.108,17.021){\makebox(0,0)[cc]{\large $\bullet$}}
\put(19.994,17.021){\makebox(0,0)[cc]{\large $\bullet$}}
\put(5.797,9.997){\makebox(0,0)[cc]{\large $\bullet$}}
\put(9.997,2.973){\makebox(0,0)[cc]{\large $\bullet$}}
\put(19.994,2.973){\makebox(0,0)[cc]{\large $\bullet$}}
\put(24.156,9.96){\makebox(0,0)[cc]{\large $\bullet$}}
\thinlines
\put(9.997,16.983){\line(3,-2){13.193}}
\put(8.176,5.983){\line(1,0){13.602}}
\put(8.176,6.02){\line(3,5){5.195}}
\put(13.371,14.679){\line(1,0){8.035}}
\put(15.051,2.973){\line(0,1){3.01}}
\put(15.051,5.983){\line(5,6){4.986}}
\put(20.037,11.966){\line(-4,5){3.954}}
\put(20.031,11.966){\line(2,-1){4.162}}
\put(13.379,14.679){\makebox(0,0)[cc]{$\bullet$}}
\put(17.875,14.642){\makebox(0,0)[cc]{$\bullet$}}
\put(16.017,16.983){\makebox(0,0)[cc]{$\bullet$}}
\put(19.994,12.004){\makebox(0,0)[cc]{$\bullet$}}
\put(19.176,10.889){\makebox(0,0)[cc]{$\bullet$}}
\put(8.176,5.946){\makebox(0,0)[cc]{$\bullet$}}
\put(15.051,5.946){\makebox(0,0)[cc]{$\bullet$}}
\put(21.777,5.983){\makebox(0,0)[cc]{$\bullet$}}
\put(21.406,14.679){\makebox(0,0)[cc]{$\bullet$}}
\put(15.088,2.936){\makebox(0,0)[cc]{$\bullet$}}
\put(23.115,8.213){\makebox(0,0)[cc]{$\bullet$}}
}
\end{picture}
\caption{Decompositions of the multiplihedron $L_3$. The left picture
shows the decomposition of $L_3$ into 3 squares corresponding to the
terms of $\p_3$. The right picture shows the decomposition of the
same multiplihedron into 10 squares  corresponding to the
terms of~$\q_3$.}
\label{fig:ob2}
\end{figure}
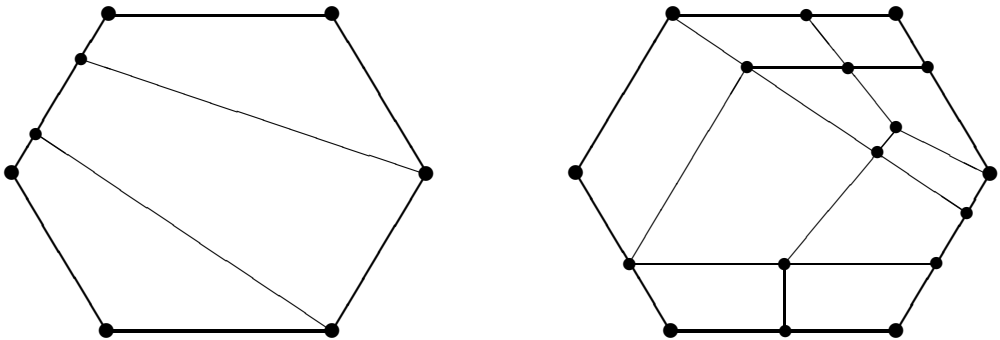

\noindent 
{\bf Transfers and polyhedra.}
The formulas for $\bfnu$, $\bfphi$, $\bfpsi$ and $\bfH$ given
in~(\ref{eq:1}) are summations of monomials in the ``initial data''
$\mu_2,\mu_3,\cdots, f,g,h$ with coefficients $\pm 1$. Ezra Getzler
conjectured that these monomials might in fact correspond to cells of
certain cell decompositions of the polyhedra governing our algebraic
structures -- Stasheff's
associahedra~\cite[page~9]{markl-shnider-stasheff:book} $K_n$, $n \geq
2$, and the multiplihedra~\cite[page~113]{markl-shnider-stasheff:book}
$L_n$, $n \geq 2$. For $K_n$, the decomposition induced by the
\pkernel\ $\p_n$ is
given by taking the tubular neighborhood of $\pa K_n$ in the
manifold-with-corners $K_n$, as illustrated for $n=4$ in
Figure~\ref{fig:ob3}.  We do not know a similar simple rule for the
multiplihedra, see also Figure~\ref{fig:ob2}.

\vskip 3mm
\noindent 
{\bf Minimal models.}
The material of this subsection is well-known to specialists.
Recall that an $A_\infty$-algebra $(W,\pa_W,\mu_2,\mu_3,\ldots)$ is
{\em minimal\/} if $\pa_W = 0$. Methods developed
in the previous sections can be 
used to construct {\em minimal models\/} of $A_\infty$-algebras as
follows.

Let $A = (V,\pa_V,\mu_2,\mu_3,\ldots)$ be an 
$A_\infty$-algebra and $W := H(V,\pa_V)$ the cohomology of its
underlying chain complex. Let $Z := {\it Ker}(\pa_V)$,
$B := {\it Im}(\pa_V)$ and 
choose a ``Hodge
decomposition'' 
\begin{equation}
\label{2}
V \cong D \oplus W \oplus B,\ \mbox { with } \ Z \cong W \oplus B.
\end{equation}
Observe that the composition 
$\omega:=  \pa_V \circ \iota_D : D \to B$, where
$\iota_D : D \hookrightarrow V$ denotes the inclusion, is a degree $-1$
isomorphism of vector spaces.
Let $f : V \to W$ be the
projection and $g: W \to V$ the inclusion induced by~(\ref{2}). 
Finally, let $h : V \to V$ be the degree $-1$ map defined as
the composition $\iota_D \circ \omega^{-1} \circ \pi_B$, where $\pi_B
: V \to B$ is the projection induced by~(\ref{2}).

It is clear that $f: (V,\pa_V) \to (W,0)$ and $g: (W,0) \to (V,\pa_V)$
are chain maps and that $h$ is a chain homotopy between $gf$ and
$\id_V$. Therefore the formula 
\begin{equation}
\label{3}
\nu_n := f \circ \p_n \circ g^{\otimes n},
\end{equation}
where
$\p_n$ is the p-kernel defined in 
Section~\ref{sec:non-induct-form}, gives a minimal model
${\mathcal M}_A  = (W,\pa_W=0,\nu_2,\nu_3,\ldots)$
of the $\ainfty$-algebra $A = (V,\pa_V,\mu_2,\mu_3,\ldots)$. 
This construction is functorial up to a choice
of the Hodge decomposition~(\ref{2}). 

More precisely, observe that decompositions~(\ref{2})
form a groupoid with morphisms given by chain
endomorphisms of $(V,\pa_V)$. We clearly have:

\begin{proposition}
The minimal model ${\mathcal M}_A$ is a functor from the groupoid
of Hodge decompositions~(\ref{2}) to the groupoid of minimal
$A_\infty$-algebras and their $A_\infty$-isomorphisms.
\end{proposition}

Observe that the ``input data'' $f,g,h$ constructed out of the Hodge
decomposition~(\ref{2}) satisfy the side conditions mentioned in
Remark~\ref{kourim_dymku}, therefore we could as well use the formulas
of~\cite{huebschmann-kadeishvili:MZ91}.

\def\cprime{$'$}

\catcode`\@=11

\vskip 1cm
\noindent
Mathematical Institute of the Academy, \v Zitn\'a 25, 115 67
Praha 1, The Czech Republic,\hfill\break\noindent
email: {\tt markl@math.cas.cz}                                                                               

\vfill
\hfill
{\tt \jobname.tex}

\end{document}